\documentclass[11pt, letterpaper]{amsart}
\usepackage{graphicx, amssymb, hyperref}

\newtheorem{thm}{Theorem}[section]

\newtheorem{prop}[thm]{Proposition}
\theoremstyle{definition}
\newtheorem{defn}[thm]{Definition}

\theoremstyle{remark}
\newtheorem{rem}[thm]{Remark}

\numberwithin{equation}{section}
\newcommand{\norm}[1]{\left\Vert#1\right\Vert}

\newcommand{\set}[1]{\left\{#1\right\}}
\newcommand{\Real}{\mathbb R}

\newcommand{\Natural}{\mathbb N}
\newcommand{\nin}{n \in \Natural}

\newcommand{\A}{\mathcal{A}}

\newcommand{\J}{\mathcal{J}}
\newcommand{\such}{{\ | \ }}
\newcommand{\limn}{\lim_{n \to \infty}}

\newcommand{\dfn}{\, := \,}

\newcommand{\prob}{\mathbb{P}}

\newcommand{\qprob}{\mathbb{Q}}
\newcommand{\expec}{\mathbb{E}}

\newcommand{\expecq}{\expec_\qprob}

\newcommand{\Lb}{\mathbb{L}}
\newcommand{\lz}{\Lb^0}
\newcommand{\li}{\Lb^\infty}
\newcommand{\lip}{\li_{+}}
\newcommand{\lzp}{\lz_{+}}

\newcommand{\F}{\mathcal{F}}

\newcommand{\ud}{\mathrm d}

\newcommand{\num}{num\'eraire}

\newcommand{\tC}{\widetilde{\C}_g}

\newcommand{\pare}[1]{\left(#1\right)}
\newcommand{\bra}[1]{\left[#1\right]}
\newcommand{\dbra}[1]{[\kern-0.15em[ #1 ]\kern-0.15em]}
\newcommand{\dbraco}[1]{[\kern-0.15em[ #1 [\kern-0.15em[}
\newcommand{\dbraoc}[1]{]\kern-0.15em] #1 ]\kern-0.15em]}
\newcommand{\C}{\mathcal{C}}

\newcommand{\Cmax}{\C^{\mathrm{max}}}
\newcommand{\K}{\mathcal{K}}

\newcommand{\Sl}{\mathcal{S}}

\newcommand{\Ll}{\mathcal{L}}

\newcommand{\Cnum}{\C^{\mathsf{num}}}

\newcommand{\indic}{\mathbb{I}}

\newcommand{\Kgc}{\mathbf{CS}_g (\C)}
\newcommand{\Koc}{\mathbf{CS}_1 (\C)}
\newcommand{\csgc}{\mathsf{cs}_g (\C)}
\newcommand{\csoc}{\mathsf{cs}_1 (\C)}
\newcommand{\csgtc}{\mathsf{cs}_1 (\tC)}

\newcommand{\absco}{{<\kern-0.53em<}}

\newcommand{\hti}{{\widetilde{h}}}

\begin{document}

\title[A structural characterization of num\'eraires of convex sets in $\lzp$]{A structural characterization of num\'eraires of convex sets of nonnegative random variables}%
\author{Constantinos Kardaras}%
\address{Constantinos Kardaras, Mathematics and Statistics Department, Boston University, 111 Cummington Street, Boston, MA 02215, USA.}%
\email{kardaras@bu.edu}%

\thanks{The author acknowledges partial support by the National Science
  Foundation, under award number DMS-0908461. Any opinions, findings and
  conclusions or recommendations expressed in this material are those
  of the author and do not necessarily reflect those of the National
  Science Foundation.}%

\date{\today}%
%\dedicatory{}%
%\commby{}%
% ----------------------------------------------------------------
\begin{abstract}
We introduce the concept of \num s of convex sets in $\lzp$, the nonnegative orthant of the topological vector space $\lz$ of all random variables built over a probability space. A necessary and sufficient condition for an element of a convex set $\C \subseteq \lzp$ to be a \num \ of $\C$ is given, inspired from ideas in financial mathematics. 
\end{abstract}

\maketitle

% ----------------------------------------------------------------

\section*{Introduction}

An element of a convex subset $\C$ in a topological vector space is called a \textsl{support point of $\C$} if it maximizes a nonzero continuous linear functional over $\C$. In finite-dimensional Euclidean spaces, every boundary point of a closed and convex set is a support point of that set. In contrast, when the topological vector space is infinite-dimensional, boundary points of a closed convex set can fail to support the set. (In fact, there exist examples of proper closed convex subsets that have \emph{no} support points --- for a specific one, see \cite{MR0149241}.)

Of immense importance, both from a probabilistic and a functional-analytic point of view, is the topological vector space $\lz$ of all (equivalence classes of real-valued) random variables built over a probability space equipped with a metric compatible with  convergence in probability. Its rich algebraic and lattice structure notwithstanding, the topological properties of $\lz$ are quite poor. In fact, if the underlying probability space is nonatomic, the topological dual of $\lz$ contains only the zero functional \cite[Theorem 2.2, page 18]{MR808777} --- in particular, convex sets in $\lz$ cannot \emph{a fortiori} have any support points according to the usual definition. In spite (and sometimes in view) of such issues, research on topological and structural properties of $\lz$ is active and ongoing; see for example \cite{MR0210177}, \cite{MR1304434}, \cite{MR1768009}, \cite{Zit08}, \cite{MR2521918}, \cite{Kar10} and \cite{KarZit10}. This note is contributing to this line of research by offering a nonstandard definition of strictly positive support points of convex sets in the nonnegative orthant of $\lz$, motivated by the well-known \num \ property in the field of financial mathematics. The main result is an interesting structural necessary and sufficient condition for a element of a convex set $\C \subseteq \lzp$ to be a \num \ of $\C$.

\section{Num\'eraires and their Structural Characterization} \label{sec: nums}

\subsection{Preliminaries}

Let $(\Omega, \F, \prob)$ be a probability space, and let $\Pi$ be
the collection of all probabilities on $(\Omega, \F)$ that are
equivalent to (the representative) $\prob \in \Pi$. Throughout the paper, $\lz$ denotes the set of all equivalence classes modulo $\Pi$ of finite real-valued random variables over $(\Omega, \F)$. We follow the usual practice of not differentiating between a random variable and the equivalence class it generates. We use $\lzp$ to denote the subset of $\lz$ consisting of elements $f \in \lz$ such that $\prob \bra{f < 0} = 0$.

The expectation of $f \in \lzp$ under $\qprob \in
\Pi$ is denoted by $\expecq [f]$. For $\qprob \in \Pi$, we define
a metric $d_\qprob$ on $\lz$ via $d_\qprob(f, g) = \expecq \bra{
  \min \set{|f - g|, 1} }$ for $f \in \lz$ and $g \in \lz$. The
topology on $\lz$ that is induced by the previous metric does not
depend on $\qprob \in \Pi$. Thus, $\lz$ becomes a complete metric space and $\lzp$ its closed subspace; convergence of sequences under the topology generated by this metric is simply convergence in $\qprob$-measure for any $\qprob \in \Pi$. Unless explicitly stated otherwise, any topological property (closedness, etc.) pertaining to subsets of $\lz$ will be understood under the aforementioned topology.

Let $\C \subseteq \lzp$. An element $f \in \C$ is called \textsl{maximal in $\C$} if the conditions $\prob[f \leq g] = 1$ and $g \in \C$ imply $\prob[f = g] = 1$; $\C^{\max}$ is then used to denote the set of all maximal elements in $\C$. Furthermore, $\C \subseteq \lzp$ will be called \textsl{bounded} if $\lim_{\ell \to \infty} \sup_{f \in \C} \qprob[f > \ell] = 0$ for some, and then for all, $\qprob \in \Pi$. The last boundedness property can be seen to coincide with boundedness of $\C$ when $\lz$ is viewed as a topological vector space \cite[Definition 5.36, page 186]{MR2378491}.

\subsection{Num\'eraires} \label{subsec: spsp}

The concept that follows is central in our development.

\begin{defn} \label{dfn: full supp point}
Let $\C \subseteq \lzp$ be convex and $g \in \C$. If $\prob \bra{f > 0, g = 0} = 0$ holds for all $f \in \C$, $g$ will be called \textsl{strictly positive on $\C$}. Furthermore, $g$ will be called a \textsl{\num \ of $\C$} if it is strictly positive on $\C$ and there exists a probability $\qprob \in \Pi$ such that $\expecq \bra{f / g \such g > 0} \leq 1$ for all $f \in \C$. The set of all \num s of $\C$ is denoted by $\Cnum$.
\end{defn}

The following result gives a more functional-analytic flavor to the concept of a \num.

\begin{prop} \label{prop: equiv num defns} 
Let $\C \subseteq \lzp$ be convex and let $g \in \C$ be strictly positive on $\C$. Then, $g$ is a \num \ of $\C$ if and only if there exists a $\sigma$-finite measure $\mu$ on $(\Omega, \, \F)$, equivalent to the probabilities in $\Pi$, such that $\int g \ud \mu = \sup_{f \in \C} \int f \ud \mu < \infty$.
\end{prop}

\begin{proof}
  We exclude from the discussion the trivial case $\C =
  \set{0}$ so that $\prob \bra{g > 0} > 0$.

First, assume that there exists a $\sigma$-finite measure $\mu$ on $(\Omega, \, \F)$, equivalent to the probabilities in $\Pi$, such that $\int g \ud \mu = \sup_{f \in \C} \int f \ud \mu < \infty$.  If $\mu[g=0]=\infty$, we can easily redefine it so that $\mu[g=0]<\infty$ without affecting the values of the integrals $\int f\ud\mu$, for $f\in\C$. Therefore, we can assume that $\mu[g=0]<\infty$. Define $\qprob \in \Pi$ via
\[
\qprob[A] = \frac{1}{2} \frac{\int_A g \ud \mu}{\int g \ud \mu} + \frac{1}{2} \frac{\mu [A \cap \{ g = 0 \}]}{\mu [\{ g = 0 \}]}, \text{ for } A \in \F,
\]
using the convention $0/0=1$. Then, $\expecq \big[ f / g \such g > 0 \big] = \int f \ud \mu / \int g \ud \mu \leq 1$ holds for all $f \in \C$.

Conversely, assume that there exists $\qprob \in \Pi$ such that $\expecq \big[ f / g \such g > 0 \big] \leq 1$ holds for all $f \in \C$. Define $\mu : \F \mapsto \Real_+ \cup \set{\infty}$ via
\[
\mu[A] = \expecq \bra{ \pare{\frac{1}{h} \indic_{\{ g > 0 \}} + \indic_{\{ g = 0 \} }} \indic_A} , \text{ for } A \in \F.
\]
It is apparent that $\mu$ is a $\sigma$-finite measure, equivalent to $\qprob \in \Pi$. Moreover, for any $f \in \C$, we have
\[
\int f \ud \mu = \expecq \big[ (f / g) \indic_{\{ g > 0 \}} \big] = \expecq [ f / g \such g > 0 ]\, \qprob[g > 0 ] \leq \qprob[g > 0] = \int g \ud \mu,
\]
which completes the proof.
\end{proof}

The previous result offers an interpretation of \num s as ``strictly positive support points'' of convex sets in $\lzp$, since $g \in \Cnum$ is supported by the ``dual'' sigma-finite measure $\mu$.  Note that the qualifying ``strictly positive'' applies both to the \num \ $g \in \C$, as well as to the supporting measure $\mu$. Of course, $g$ is not a support point of $\C$ in the traditional functional-analytic sense, since the mapping $\lzp \in f \mapsto \int f \ud \mu$ is only lower semi-continuous. However, when $\C$ is viewed as a convex set in the Banach space $\Lb^1(\Omega, \F, \mu)$, then $g$ \emph{is} a support point of $\C$ in the usual sense.

\smallskip

Let $\C \subseteq \lzp$ be convex, and let $g \in \C$ be strictly positive on $\C$. The question we focus on is the following: \emph{Is there a structural equivalent to the condition that $g$ is a \num \ of $\C$?} Necessary conditions are easy to obtain. For example, if $\C$ is to afford \emph{any} \num s, then $\C$ has to be bounded, as it immediately follows by a use of Chebyshev's inequality. Also, if $g \in \Cnum$, it is clearly necessary that $g \in \Cmax$. As we shall shortly see in \S \ref{subsec: exa}, the previous two necessary conditions ($\C$ is bounded and $g \in \Cmax$) are not always sufficient to ensure that $g \in \Cnum$. A detailed understanding of the issues faced in an example presented in \S \ref{subsec: exa} below will enable us to eventually reach our main result, Theorem \ref{thm: main}.

\subsection{An example} \label{subsec: exa}

In financial mathematics, a convex set $\C \subseteq \lzp$ consisting of terminal values of nonnegative stochastic integrals starting from unit initial value with respect to a semimartingale integrator is used to model discounted outcomes of wealth processes starting from unit capital. More precisely, the semimartingale integrator models discounted asset prices and the predictable integrands model investment strategies. As long as there are no constrains on investment (further from the natural constraints of nonnegativity for the involved wealth processes), and when $g \in \C$ is such that $\prob \bra{g > 0} = 1$, the condition that $\C$ is bounded and $g \in \Cmax$ is, quite interestingly, \emph{equivalent} to $g \in \Cnum$. (See \cite{MR1304434, MR1381678} for a comprehensive treatment of this topic.) However, in the presence of investment constraints, the situation becomes more complicated, as we present below with an illustrating example.

\smallskip

Start with a probability space $(\Omega, \F, \prob)$, rich enough to support $\xi \in \lzp$ with $\prob[\xi > 0] = 1$, and $\prob [\xi \leq \epsilon] > 0$ as well as $\prob [1 / \xi \leq \epsilon] > 0$ holding for all $\epsilon > 0$. Let $S_i = (S_i(t))_{t \in \set{0,T}}$ for $i \in \set{1, 2}$ be defined via $S_1(0) = 1 = S_2(0)$, and $S_1(T) = \xi$, $S_2(T) = 1 + \xi$. Each $S_i$, $i \in \set{1, 2}$ is modeling the discounted price of a financial asset. For any $\vartheta = (\vartheta_1, \vartheta_2) \in \Real^2$, define $X^{\vartheta}$ via $X^\vartheta (0) = 1$ and
\[
X^{\vartheta}(T) = 1 + \vartheta_1 (S_1(T) - S_1(0)) + \vartheta_2 (S_2(T) - S_1(0)) = 1 - \vartheta_1 + (\vartheta_1 + \vartheta_2) \xi;
\]
then, $X^{\vartheta}(T)$ is modeling the discounted financial outcome at time $T$ of an investment starting with unit capital and holding a position $\vartheta$ in the assets.

We now introduce constraints on investment. Let $C \dfn \set{(\vartheta_1, \vartheta_2) \in \Real_+^2 \such \vartheta_2 \leq \sqrt{\vartheta_1} \leq 1}$, which is a convex and compact subset of $\Real_+^2$. It is easy to check that $X^{\vartheta}(T) \geq 0$, for all $\vartheta \in C$. Consider
\[
\C = \set{X^{\vartheta}_T  \such \vartheta \in C} = \set{1 - \vartheta_1 + (\vartheta_1 + \vartheta_2) \xi \such \vartheta \in C},
\]
which is a convex, closed and bounded subset of $\lzp$. Using the fact that $\prob [\xi \leq \epsilon] > 0$ and $\prob [1 / \xi \leq \epsilon] > 0$ hold for all $\epsilon > 0$, it is straightforward to check that
\begin{equation} \label{eq: C_max_exa}
\Cmax = \set{1 - \gamma + (\gamma + \sqrt{\gamma}) \xi \such \gamma \in [0, 1]};
\end{equation}
in particular, $1 \in \Cmax \subseteq \C$.

Although both $\C$ is bounded and $1 \in \Cmax$ hold, we claim that $1 \notin \Cnum$. To this end, suppose that $\qprob \in \Pi$ is such that $\expecq \bra{f} \leq 1$ for all $f \in \C$. Then, $\expecq \bra{1 - \gamma + (\gamma + \sqrt{\gamma}) \xi} \leq 1$ for all $\gamma \in [0, 1]$. Rearranging, $\expecq \bra{\xi} \leq  \gamma / \pare{\gamma + \sqrt{\gamma}} = \sqrt{\gamma} / \pare{\sqrt{\gamma} + 1}$, for all $\gamma \in \, ]0, 1]$. This would imply that $\expecq \bra{\xi} = 0$, i.e., $\qprob \bra{\xi > 0} = 0$ which, in view of $\prob[\xi > 0] = 1$, clearly contradicts the equivalence between $\prob$ and $\qprob$.

\begin{rem} \label{rem: key}
It is worthwhile to try to understand what structural property of $\C$ prevented $g=1$ from being a \num \ of $\C$ in the above example, since it will help shed light on the exact necessary and sufficient conditions needed in the statement of our main result. For the time being, consider any $\C \subseteq \lzp$ and any $g \in \C$. Suppose that $g \in \Cnum$, and let $\qprob \in \Pi$ be as in Definition \ref{dfn: full supp point}. Pick $f \in \C$ and $\delta \in \Real_+$ such that, with $f' \dfn (1 + \delta) f - \delta g$, we have $f' \in \lzp$. If $\C$ represents terminal outcomes from investment as in the example above, $f'$ corresponds to taking a long position of $(1 + \delta)$ units of the portfolio leading to the outcome $f$ and a short position on $\delta$ units of the portfolio leading to the outcome $g$; the fact that $f' \in \lzp$ guarantees that there is no risk of going negative. As $\expecq \bra{f / g \such g > 0} \leq 1$, we obtain $\expecq \bra{f' / g \such g > 0} \leq 1$ as well. Note that the previous holds for all possible $f' \in \lzp$ constructed as before. The upshot is the following: if we enlarge $\C$ by including all such combinations (taking short positions on $g$), a use of Chebyshev's inequality implies that we still end up with a set that is bounded. The previous observation, however simple, will be key in the development.
\end{rem}

We return to our concrete example. For $\nin$, let $f_n \dfn n / (1+n) + \pare{1 / (1+n) + 1 / \sqrt{1+n}} \xi$; by \eqref{eq: C_max_exa}, $f_n \in \Cmax \subseteq \C$. With $f'_n \dfn (1+n) f_n - n = (1 + \sqrt{1+n}) \xi$, we have $f'_n \in \lzp$ for all $\nin$. By the discussion in Remark \ref{rem: key} above, if $1$ were to be a \num \  of $\C$, $\set{f'_n \such \nin}$ would have to be a bounded subset of $\lzp$, which is plainly false.

\subsection{The equivalence result} \label{subsec: enlargement}

Guided by the discussion of \S \ref{subsec: exa}, for $\C \subseteq \lzp$ and $g \in \C$ we define $\Kgc$ as the class of all $\K \subseteq \lzp$ such that:
\begin{enumerate}
\item[(CS1)] $\C \subseteq \K$.
\item[(CS2)] $\K$ is convex and closed.
\item[(CS3)] If $f \in \K$ and $\delta \in \Real_+$ are such that $\pare{(1 + \delta) f - \delta g} \in \lzp$, then $\pare{(1 + \delta) f - \delta g} \in \K$.
\end{enumerate} 
It is clear that $\Kgc$ is closed under arbitrary intersections. Furthermore, $\lzp \in \Kgc$, i.e., $\Kgc \neq \emptyset$. Therefore, there exists a minimal set in $\Kgc$, which we shall denote by $\csgc$:
\begin{equation} \label{eq: csgc}
\csgc \dfn \bigcap_{\K \in \Kgc} \K.
\end{equation}

The combination of (CS1) and (CS2) plainly states that sets in $\Kgc$ are closed and convex enlargements of $\C$. Using jargon from financial mathematics,  (CS3) states that these enlargements of $\C$ are at least large enough to contain all results from leveraged positions using short selling of $g$, so long as these combinations lead to nonnegative outcomes. The minimal way of doing so is given by the set $\csgc$ of \eqref{eq: csgc}. (In $\csgc$, ``$\mathsf{c}$'' is used as a mnemonic for \emph{closed and convex} and ``$\mathsf{s}_g$'' as a mnemonic for \emph{short sales in $g$}.) As there does not seem to exist a constructive way to obtain $\csgc$ from $\C$,  \eqref{eq: csgc} is utilized as its definition.

\smallskip

After all the preparation, we are ready to state our main equivalence result.

\begin{thm} \label{thm: main}
Let $\C \subseteq \lzp$ be convex, and let $g \in \C$ be strictly positive on $\C$. Define $\csgc$ as in \eqref{eq: csgc}. Then, $g \in \Cnum$ if and only if $\csgc$ is bounded.
\end{thm}

If $g \in \Cnum$, it easily follows that $\csgc$ is bounded. Indeed, pick a $\qprob \in \Pi$, such that $\expecq \bra{f / g \such g > 0} \leq 1$ for all $f \in \C$. Define
\[
\K \dfn \set{h \in \lzp \such \prob \bra{h > 0, \, g = 0} = 0 \text{ and } \expecq \bra{h / g \such g > 0} \leq 1}.
\]
It is straightforward to check that $\K \in \Kgc$; this means that $\csgc \subseteq \K$. By Chebyshev's inequality, $\K$ is bounded; therefore, $\csgc$ is bounded as well. The proof of the more involved converse implication is discussed in Section \ref{sec: proof} below.

\section{The Proof of Theorem \ref{thm: main}} \label{sec: proof}

Let $\C \subseteq \lzp$ be convex, and let $g \in \C$ be strictly positive on $\C$. Assume that  $\csgc$ is bounded. To complete the proof of Theorem \ref{thm: main}, we have to establish that $g \in \Cnum$. In order to ease the reading and understanding, we split the proof in four steps.

\subsection*{Step 1}

We begin by showing the \emph{we can reduce the proof to the case $g=1$}. Indeed, define the convex set $\tC \dfn \set{\indic_{\{ g = 0 \}} + \indic_{\{ g > 0 \}} (f / g) \such f \in \C}$. Then, $1 \in \tC$ is strictly positive on $\tC$. One can check that $\csgtc = \set{\indic_{\{ g = 0 \}} + \indic_{\{ g > 0 \}} (h / g) \such h \in \csgc}$. This implies that $\csgc$ is bounded if and only if $\csgtc$ is bounded. Now, suppose that $1 \in \tC$ is a \num \ of $\tC$; in other words, that there exists $\qprob \in \Pi$ such that $\expecq \bra{f} \leq 1$ holds for all $f \in \tC$. The last is equivalent to $\expecq \bra{ (f / g) \indic_{\set{g > 0}}  }  \leq \qprob[g > 0]$ holding for all $f \in \C$, which shows that $g \in \Cnum$. 

In view of the above discussion, \emph{we assume from now on until the end of the proof that $g=1$}.

\subsection*{Step 2}
Define $\Sl \dfn \set{f \in \lzp \such 0 \leq f \leq h \text{ for some } h \in \csoc}$ be the \emph{solid hull} of $\csoc$.  We shall show below that \emph{$1 \in \Sl^{\max}$ and $\Sl \in \Koc$}.

Clearly, $1 \in \Sl^{\max}$ is equivalent to $1 \in \csoc^{\max}$. Suppose then that $f \in \csoc$ is such that $\prob \bra{f \geq 1} = 1$. By property (CS3) of the sets in $\Koc$ mentioned in \S \ref{subsec: enlargement}, we have $f_n \dfn f + n(f-1) = \pare{(n+1) f - n} \in \csoc$ for all $\nin$. If $\prob[f > 1] > 0$, the $\csoc$-valued sequence $(f_n)_{\nin}$ would fail to be bounded. Therefore, $\prob [f = 1] = 1$, which implies that $1 \in \csoc^{\max}$.

We proceed in showing that $\Sl \in \Koc$. We have $\C \subseteq \csoc \subseteq \Sl$, which shows that $\Sl$ satisfies property (CS1). Further, it is straightforward to check that $\Sl$ is convex and bounded. It is also true that $\Sl$ is closed. (To see the last fact, pick an $\Sl$-valued sequence  $(f_n)_{n \in \Natural}$ that converges to $f \in \lz_{+}$; we need to show that $f \in \Sl$. By passing to a subsequence if necessary, we may assume that $\prob \bra{\limn f_n = f} = 1$. Let $(\hti_n)_{n \in \Natural}$ be a $\csoc$-valued sequence with $\prob \big[ f_n \leq \hti_n \big] = 1$ for all $n \in \Natural$. By \cite[Lemma A1.1]{MR1304434}, we can extract a sequence $(h_n)_{n \in \Natural}$ such that, for each $n \in \Natural$, $h_n$ is a convex combination of $\hti_n, \hti_{n+1}, \ldots$, as well as $\prob \bra{\limn h_n = h} = 1$ holds for some $h \in \lzp$. Of course, $h \in \csoc$ and it is easy to see that $\prob \bra{f \leq h} = 1$. We then conclude that $f \in \Sl$.) This shows that $\Sl$ satisfies property (CS2). Now, let $f \in \Sl$ be such that $\pare{(1 + \delta) f - \delta } \in \lzp$ for some $\delta \in \Real_+$. Pick $h \in \csoc$ with $\prob \bra{f \leq h} = 1$. Then, $\pare{(1 + \delta) h - \delta } \in \lzp$ also holds. By definition of $\csoc$, we have $\pare{(1 + \delta) h - \delta } \in \csoc$. As $\pare{(1 + \delta) f - \delta g} \in \lzp$ and $\prob \bra{(1 + \delta) f - \delta  \leq (1 + \delta)h - \delta } = 1$, we obtain that $\pare{(1 + \delta) f - \delta } \in \Sl$; therefore, $\Sl$ also satisfies property (CS3). We conclude that $\Sl \in \Koc$.

\subsection*{Step 3}

In the sequel, $\li$ denotes the space of essentially bounded (modulo $\prob$) elements of $\lz$. Note that topological notions are \emph{still} considered under $\lz$.

Define $\Ll \dfn \Sl \cap \li$. All the statements regarding $\Ll$ below, which we shall be using tacitly, follow in a straightforward way from the properties of $\Sl$:
\begin{itemize}
\item $\Ll$ is convex and solid. (The latter means that $0 \leq f \leq g\in \Ll$ implies $f \in \Ll$.)
\item $1 \in \Ll^{\max}$.
\item For all $f \in \Sl$ and $\nin$, $\min \set{f,n} \in \Ll$.
\item For any uniformly bounded (modulo $\prob$) $\Ll$-valued sequence $(f_n)_{\nin}$ that converges to $f \in \lzp$, we have $f \in \Ll$.
\item If $f \in \Ll$ and $\delta \in \Real_+$ are such that $\pare{(1 + \delta) f - \delta} \in \lzp$, then $\pare{(1 + \delta) f - \delta} \in \Ll$.
\end{itemize}
Continuing, define $\A_\alpha \dfn \alpha (\Ll - 1) = \set{\alpha(f - 1) \such f \in \Ll}$ for $\alpha \in \Real_+$, as well as $\J \dfn \bigcup_{\alpha \in \Real_+} \A_\alpha$. We shall show below that \emph{$\J$ is a weak*-closed convex cone in $\li$, satisfying $\J = \J - \lip$ and $\J \cap \lip = \set{0}$}. (We obviously define $\lip \dfn \lzp \cap \li$; furthermore, the weak* topology on the Banach space $\li$ equipped with the usual $\li$-norm is defined as usual.)

It is clear that $\J$ is a convex cone in $\li$. Also, since $1 \in \Ll^{\max}$, $\J \cap \lip = \set{0}$ is immediate. We proceed in showing that $\J = \J - \lip = \set{\phi - h \such \phi \in \J, \, h \in \lip}$. Since $\J \subseteq \J - \lip$, we only have to show that if $\psi = \phi - h$ where $\phi \in \J$ and $h \in \lip$, then $\psi \in \J$. We assume that $\prob[ h > 0] > 0$; otherwise, $\psi \in \J$ is trivial. Write $\phi = \alpha (f - 1)$, where $f \in \Ll$ and $\alpha \in \Real_+$.  With $\eta = \norm{h}_{\infty} \in \Real_+$, so that $\prob[h \leq \eta] = 1$, let $f' \dfn \eta / (\alpha + \eta) + \pare{\alpha / (\alpha + \eta)} f$. Since $1 \in \Ll$, $f \in \Ll$, and $\Ll$ is convex, we have $f' \in \Ll$. Now, define $f'' \dfn (\eta - h) / (\alpha + \eta) + \pare{\alpha / (\alpha + \eta)} f$; then $f'' \in \lip$ and $f'' \leq f'$; since $\Ll$ is solid, $f'' \in \Ll$. Then,
\[
\psi = \phi - h = \alpha(f-1) - h = (\alpha + \eta) (f'' - 1) \in \A_{\alpha + \eta} \subseteq \J,
\]
which establishes our claim $\J = \J - \lip$.

It only remains to establish that $\J$ is weak*-closed in $\li$. Before this is done, we show that $\psi \in \J$ and $\prob \bra{\psi \geq -1} = 1$ imply $\psi \in \A_1$. First, note that $\A_\alpha \subseteq \A_\beta$ whenever $0 \leq \alpha < \beta$: indeed, for $f \in \Ll$, use the fact that $\Ll$ is convex to write $\alpha (f - 1) = \beta (f' - 1)$, where
\[
f' \dfn \pare{\frac{\alpha}{\beta} f + \frac{\beta - \alpha}{\beta}} \in \Ll.
\]
Now, let $\psi \in \A_\alpha$ with $\prob \bra{\psi \geq -1} = 1$. If $\alpha \leq 1$, $\psi \in \A_1$ is obvious by the above discussion. Assume that $\alpha > 1$. Write $\psi = \alpha(f - 1)$ for $f \in \Ll$ and note that $\prob \bra{\psi \geq -1} = 1$ translates to $\pare{\alpha f - (\alpha - 1)} \in \lzp$. In that case, with $f' \dfn \alpha f - (\alpha - 1)$ we have $f' \in \Ll$. But then, $\psi = \alpha (f-1) = (f' - 1)  \in \A_1$.

We shall now show that $\J$ is weak*-closed in $\li$. By combining the Krein-Smulian theorem with the fact that, for uniformly bounded and convex sets of $\li$, weak*-closedness coincides with $\lz$-closedness (in this respect, see also \cite[Theorem 2.1]{MR1304434}), it suffices to show that for any $\J$-valued sequence $(\phi_n)_{\nin}$ that converges (in $\lz$) to $\phi \in \lz$ and is such that $\prob \bra{|\phi_n| \leq 1} = 1$ for all $\nin$, we have $\phi \in \J$. As $\prob \bra{\phi_n \geq - 1} = 1$ for all $\nin$, $(\phi_n)_{\nin}$ is $\A_1$-valued by the discussion of the preceding paragraph. For $\nin$, write $f_n = \phi_n + 1$; then, $(f_n)_{\nin}$ is $\Ll$-valued, it converges to $f \dfn \phi + 1$ and $\prob \bra{0 \leq f_n \leq 2} = 1$ for all $\nin$. From the properties of $\Ll$, it follows that $f \in \Ll$, i.e., that $\phi \in \A_1 \subseteq \J$.

\subsection*{Step 4}
We can now finish the proof of Theorem \ref{thm: main}. Using Step 3 above, an invocation of the Kreps-Yan separation theorem (see \cite{MR611252} and \cite{MR580127}) gives the existence of a probability $\qprob \in \Pi$, such that $\expecq \bra{\phi} \leq 0$ holds for all $\phi \in \J$. (Expectations under probabilities in $\Pi$ of elements of $\li$ are always well-defined.) It follows that $\expecq \bra{f} \leq 1$ for all $f \in \Ll$. For $f \in \Sl$, we have $\min \set{f, n} \in \Ll$ for all $\nin$. Then, $\expecq \bra{f} = \limn \expecq \bra{\min \set{f, n}} \leq 1$ holds for all $f \in \Sl$; in other words, $\expecq \bra{f} \leq 1$ for all $f \in \Sl$. Finally, since $\C \subseteq \csoc \subseteq \Sl$, we obtain $\expecq \bra{f} \leq 1$ for all $f \in \C$.

%----------------------------------------------------------------
\bibliographystyle{siam}
\bibliography{struct_num}
\end{document}